\documentclass[12pt]{amsart}

\usepackage{amsfonts, amsmath, amssymb,amsthm,mathrsfs}
\usepackage[dvips]{graphicx, color}

\newcommand{\vanish}[1]{}

\newcommand{\Z}{\mathbb{Z}}  
 
\newcommand{\QSym}{\mathcal{Q}}
 
\newcommand{\rk}{\rho}  
\newcommand{\F}{\mathcal{F}} 
 
\newcommand{\ext}{\operatorname{ext}}  
\newcommand{\closure}[1]{\left\langle #1 \right\rangle}
\newcommand{\lk}{\mathit{lk}}
 
\newtheorem{theorem}{Theorem}
\newtheorem{proposition}[theorem]{Proposition}
\newtheorem{lemma}[theorem]{Lemma}
\newtheorem{corollary}[theorem]{Corollary}

\theoremstyle{definition}
\newtheorem{definition}[theorem]{Definition}

\theoremstyle{remark}
\newtheorem{remark}[theorem]{Remark} 

\numberwithin{theorem}{section}
\numberwithin{equation}{section}

\begin{document}

\title[Convex geometries and their spheres]{Enumeration in convex geometries and associated polytopal subdivisions
of spheres}
 
\author{Louis J.\ Billera}
\address{Department of Mathematics, Cornell University, Ithaca, NY
14853-4201}
\email{billera@math.cornell.edu}
\thanks{The first author was supported in part by NSF grant
DMS-0100323.  The second author was supported by an NSF
Postdoctoral Fellowship.  The first two authors enjoyed the hospitality of
the Mittag-Leffler Institute, Djursholm, Sweden, during the preparation
of this manuscript.}

\author{Samuel K. Hsiao}
\address{Department of Mathematics, University of Michigan, Ann Arbor,
 MI 48109-1043}
 \curraddr{Department of Mathematics, Bard College, P.O. Box 5000, Annandale-on-Hudson, NY 12504-5000}
\email{shsiao@umich.edu}

\author{J. Scott Provan}
\address{Department of Statistics and Operations Research, University of North Carolina,
Chapel Hill, NC 27599-3180}
\email{Scott\_Provan@UNC.edu}


\begin{abstract}
We construct  CW spheres from the lattices that arise as the closed sets of a convex
closure, the meet-distributive lattices.
These spheres are nearly polytopal, in the sense that their barycentric
subdivisions are simplicial polytopes.
The complete information on the numbers of faces and chains of faces
in these spheres can be obtained from the defining lattices in a manner
analogous to the relation between arrangements of hyperplanes and
their underlying geometric intersection lattices.
\end{abstract}

\maketitle


\section{Introduction}
A well known result due to Zaslavsky \cite{Z} shows that the numbers of faces
in an arrangement of hyperplanes in a real Euclidean space can be read
from the underlying geometric lattice of all intersections of these
hyperplanes.  This result was extended to the determination of the numbers of
chains of faces in arrangements in \cite{BS,matroid}.  In \cite{BER}, the numbers
of chains in an arrangement were shown to depend only on the numbers of chains
in the associated geometric lattice.  A particularly simple form of this relationship, in terms of
quasisymmetric functions,  was given in \cite{BHW}.

Geometric lattices (matroids) are combinatorial abstractions of linear span in vector
spaces.  There is a different combinatorial model for convex span, known as convex
geometries (or anti-matroids)  \cite{Edel1,EdelJam,Edel2,EdelSaks}, for which the
corresponding lattices are the meet-distributive lattices.  We show here
that a similar situation exists for these; that is, for each convex geometry, we construct a
regular CW sphere, whose enumerative properties are related to those of the underlying
geometry in essentially the same way.  Moreover, these spheres are nearly polytopes,
in the sense that their first barycentric subdivisions are combinatorially simplicial
convex polytopes.

We begin by establishing some notation.  Our basic object of study is
a combinatorial closure operation called a {\em convex} or
{\em anti-exchange} closure.  This is defined on a finite set, which
we will take, without loss of generality, to be the set
$[n]:= \{ 1,2,\dots ,n \}$.  
\begin{definition}
A {\em convex closure} is a
function $\langle \cdot \rangle : 2^{[n]} \rightarrow 2^{[n]}$,
 $A \mapsto \langle A \rangle$,
such that, for $A,B \subseteq [n]$,
\begin{enumerate}
\item
$A \subseteq \langle A \rangle$
\item
if $A \subseteq B$ then $\langle A \rangle \subseteq \langle B \rangle$
\item
$\langle A \rangle = \langle \langle A \rangle \rangle $
\item
if $x,y \notin \langle A \rangle$ and $x \in \langle A \cup {y} \rangle$
then $y \notin \langle A\cup {x} \rangle$.
\end{enumerate}
\end {definition}

The last condition is often called the {\em anti-exchange} axiom, and
the complements of the closed sets of such
a closure system has been called an {\em anti-matroid}.  We will call a
set together with a convex closure operator on it a {\em convex geometry}.
The set of {\em closed sets} of a convex geometry, that is, those sets $A$
satisfying $A=\langle A \rangle$, form a lattice when ordered by set inclusion.
Such lattices are precisely the  meet-distributive lattices.  A lattice
$L$ is {\em meet-distributive} if for each $y \in L$, if $x \in L$ is
the meet of
(all the) elements covered by $y$, then the interval $[x,y]$ is a Boolean
algebra.

\begin{figure}[htb]
\begin{center}
\setlength{\unitlength}{2000sp}%
\begingroup\makeatletter\ifx\SetFigFont\undefined%
\gdef\SetFigFont#1#2#3#4#5{%
  \reset@font\fontsize{#1}{#2pt}%
  \fontfamily{#3}\fontseries{#4}\fontshape{#5}%
  \selectfont}%
\fi\endgroup%
\begin{picture}(6955,3516)(-3685,-3061)
\thicklines
\put(-1200,-1580){\makebox(0,0)[lb]{\smash{{\SetFigFont{12}{20.4}{\rmdefault}{\mddefault}{\itdefault}{\color[rgb]{0,0,0}$c$}%
}}}}
{\color[rgb]{0,0,0}\put(600,-660){\circle{80}}
}%
{\color[rgb]{0,0,0}\put(2400,-660){\circle{80}}
}%
{\color[rgb]{0,0,0}\put(600,-1560){\circle{80}}
}%
{\color[rgb]{0,0,0}\put(1500,-2460){\circle{80}}
}%
{\color[rgb]{0,0,0}\put(2400,-1560){\circle{80}}
}%
{\color[rgb]{0,0,0}\put(1500,-1560){\circle{80}}
}%
{\color[rgb]{0,0,0}\put(-3600,-1260){\circle{80}}
}%
{\color[rgb]{0,0,0}\put(-2400,-1260){\circle{80}}
}%
{\color[rgb]{0,0,0}\put(-1200,-1260){\circle{80}}
}%
{\color[rgb]{0,0,0}\put(1500,240){\line(-1,-1){900}}
\put(600,-660){\line( 0,-1){900}}
\put(600,-1560){\line( 1,-1){900}}
\put(1500,-2460){\line( 1, 1){900}}
\put(2400,-1560){\line( 0, 1){900}}
\put(2400,-660){\line(-1, 1){900}}
}%
{\color[rgb]{0,0,0}\put(600,-660){\line( 1,-1){900}}
\put(1500,-1560){\line( 1, 1){900}}
}%
{\color[rgb]{0,0,0}\put(1500,-1560){\line( 0,-1){900}}
}%
{\color[rgb]{0,0,0}\put(-3600,-1260){\line( 1, 0){2400}}
}%
\put(1500,-3060){\makebox(0,0)[lb]{\smash{{\SetFigFont{14}{34.8}{\rmdefault}{\mddefault}{\itdefault}{\color[rgb]{0,0,0}$L$}%
}}}}
\put(1801,239){\makebox(0,0)[lb]{\smash{{\SetFigFont{12}{24.0}{\rmdefault}{\mddefault}{\itdefault}{\color[rgb]{0,0,0}$\{a,b,c\}$}%
}}}}
\put(901,-661){\makebox(0,0)[lb]{\smash{{\SetFigFont{12}{24.0}{\rmdefault}{\mddefault}{\itdefault}{\color[rgb]{0,0,0}$\{a,b\}$}%
}}}}
\put(2600,-661){\makebox(0,0)[lb]{\smash{{\SetFigFont{12}{24.0}{\rmdefault}{\mddefault}{\itdefault}{\color[rgb]{0,0,0}$\{b,c\}$}%
}}}}
\put(650,-1561){\makebox(0,0)[lb]{\smash{{\SetFigFont{12}{24.0}{\rmdefault}{\mddefault}{\itdefault}{\color[rgb]{0,0,0}$\{a\}$}%
}}}}
\put(1651,-1636){\makebox(0,0)[lb]{\smash{{\SetFigFont{12}{24.0}{\rmdefault}{\mddefault}{\itdefault}{\color[rgb]{0,0,0}$\{b\}$}%
}}}}
\put(2600,-1561){\makebox(0,0)[lb]{\smash{{\SetFigFont{12}{24.0}{\rmdefault}{\mddefault}{\itdefault}{\color[rgb]{0,0,0}$\{c\}$}%
}}}}
\put(1801,-2536){\makebox(0,0)[lb]{\smash{{\SetFigFont{12}{24.0}{\familydefault}{\mddefault}{\updefault}{\color[rgb]{0,0,0}$\emptyset$}%
}}}}
\put(-3600,-1580){\makebox(0,0)[lb]{\smash{{\SetFigFont{12}{20.4}{\rmdefault}{\mddefault}{\itdefault}{\color[rgb]{0,0,0}$a$}%
}}}}
\put(-2350,-1580){\makebox(0,0)[lb]{\smash{{\SetFigFont{12}{20.4}{\rmdefault}{\mddefault}{\itdefault}{\color[rgb]{0,0,0}$b$}%
}}}}
{\color[rgb]{0,0,0}\thicklines
\put(1501,239){\circle{80}}
}%
\end{picture}%
\end{center}
\caption{The convex geometry of three collinear points and its associated meet-distributive
lattice $L$.\label{threepoints}}
\end{figure}
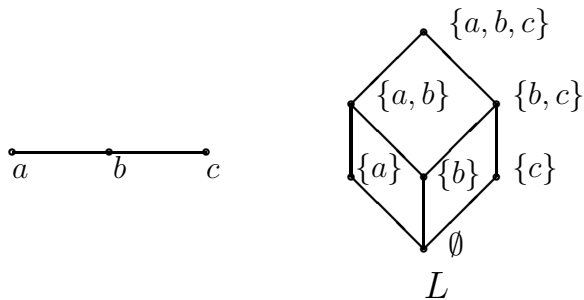

One example of an anti-exchange closure operator is ideal closure on a
partially ordered set $P$; here, for $A\subseteq P$, $\langle A \rangle$
denotes the (lower) order ideal generated by $A$.  The lattices of closed
sets of these are precisely the distributive lattices.  Another class of examples
comes from considering convex closure on a finite point set in Euclidean space.
Figure \ref{threepoints} illustrates the convex geometry formed by three collinear
points $a,b,c$.  Note that the set $\{a,c\}$ is not closed since its closure is $\{a,b,c\}$.

Meet-distributive lattices were first studied by Dilworth \cite{Dil} and
have reappeared in many contexts since then (see \cite{Monj}).
Their study in the context of theory of convex geometries was
extensively developed about 20 years ago in a series of papers
 by Edelman and coauthors \cite{Edel1,EdelJam,Edel2,EdelSaks}.  See also
\cite{KL} and \cite{BZ} for general discussions.  An important (and
characterizing) property of convex geometries is that every set has
a unique minimal generating set, that is, for each $A \subseteq [n]$,
there is a unique minimal subset $\ext(A) \subseteq A$ so that 
$\langle A \rangle = \langle \ext(A) \rangle$ \cite[Theorem 2.1]{EdelJam}.
The elements of $\ext(A)$ are called the {\em extreme points} of $A$.

A {\em simplicial complex} on a finite set $V$ is a family of subsets
$\Delta \subseteq 2^V$ such that if $\tau \subseteq \sigma \in \Delta$
then $\tau \in \Delta$ and $\{v\}\in \Delta$ for all $v\in V$.  The
elements of $\Delta$ are called the {\em faces} of the complex and the
elements of $V$ are its {\em vertices}. 
We will need an operation on simplicial complexes known as
{\em stellar subdivision}.

\begin{definition}
The {\em stellar subdivision} of a simplicial complex $\Delta$ over
a nonempty face $\sigma \in \Delta$ is the simplicial complex 
$sd_\sigma(\Delta)$ on the set $V\cup \{v_\sigma\}$, where $v_\sigma$
is a new vertex, consisting of
\begin{enumerate}
\item
all $\tau \in \Delta$ such that $\tau \not\supseteq \sigma$, and
\item
all $\tau \cup \{v_\sigma\}$ where $\tau \in \Delta$,
 $\tau \not\supseteq \sigma$ and $\tau \cup \sigma \in \Delta$.
\end{enumerate}
\end{definition}

The use of stellar subdivision to describe order complexes of
posets was begun in \cite{Pr}, where it  was shown that the order
complex of any distributive lattice can be obtained from a simplex
by a sequence of stellar subdivisions.  Although this result and some
of its implications were discussed in \cite{PrB}, its proof was never
published.  We will give a generalization of this result to
meet-distributive lattices in the next section.  The proof is an
adaptation of that in \cite{Pr}.

More recently, stellar subdivision was used in \cite{BPSZ} to produce the
order complex of a so-called Bier poset of a poset $P$ from the
order complex of $P$.

In \S 2, we describe the order complex of a meet-distributive lattice
as a stellar subdivision of a simplex.  We use this in \S 3 to construct the
sphere associated with the lattice.  Finally, in \S 4 we relate enumeration
in this sphere to that of the lattice.


\section{Order complexes of meet-distributive lattices}
Let $L$ be an arbitrary meet-distributive lattice.  We can assume $L$ is
the lattice of closed sets of a convex closure $\closure{\cdot}$ on
the set $[n]$, for some $n>0$.  $L$ has unique maximal element
$\hat{1} = \closure{[n]}$ and minimal element
$\hat{0} = \closure{\emptyset}$ (we may assume
$\closure{\emptyset}=\emptyset$, although this will not be important
here).  For simplicity of notation, we will write $\closure{i}$ for
the {\em principal} closed set $\closure{\{i\}}$ whenever $i\in [n]$.
These are precisely the {\em join-irreducible} elements of $L$, that is,
those $x\in L\setminus \{\hat{0}\}$ that cannot be written as 
$y \vee z$, with $y,z<x$.  (This follows, for example,
from \cite[Theorem 2.1(f)]{EdelJam}.)

In fact, the convex closure
$\closure{\cdot}$ is uniquely defined from the lattice $L$: we take
$[n]$ to be an enumeration of the join-irreducible elements of $L$
and define, for $A\subseteq [n]$, 
$$\closure{A} = \left\{ j\in [n] ~\Big\vert ~ j \le \bigvee_{i\in A}i \right\}.$$
Thus we are free, without loss of generality, to use the closure relation
when making constructions concerning the lattice $L$.

Consider the simplex of all principal closed sets (join-irreducibles)
  $\{ \closure{i}~ |~ i \in [n] \}$,
and let $\Delta_0$ be the simplicial complex consisting of this simplex and all
its faces (subsets).
Note that $\Delta(L\setminus \{\hat{0}\})$, the order complex of
$L\setminus \{\hat{0}\}$, is a simplicial
complex on the vertex set
$V=\{ \closure{A} | A\subseteq [n], A\ne \emptyset \}$.

\begin{theorem}
For any meet-distributive lattice $L$,  $\Delta(L\setminus \{\hat{0}\})$
can be obtained from the simplex of join-irreducible elements by a sequence
of stellar subdivisions.
\label{ordercomplex}
\end{theorem}
\begin{proof}
Suppose $L$ is the lattice of closed sets of a convex closure $\closure{\cdot}$ on $[n]$.
Let $A_{1}, A_{2}, \dots , A_{k}$ be a reverse linear extension of $L\setminus \{\hat{0}\}$, that is,
the $A_{i}$ are all the nonempty closed subsets in $[n]$, ordered so that we never
have $A_{i} \subseteq A_{j}$ if $i<j$.  In particular, $A_{1}=[n]$.

The order complex $\Delta(L\setminus \{\hat{0}\})$ can be obtained from $\Delta_0$
by a sequence of stellar subdivisions as follows.  For $i=1,\dots,k$, let
$$\Delta_{i} = sd_{\ext(A_{i})}(\Delta_{i-1}),$$
where, by a slight abuse of notation, $\ext(A_{i})$ will denote
the face of $\Delta_{i-1}$ having vertices $\closure{j}$,
$j \in \ext(A_{i})$.  The new vertex added at  the $i^{th}$ step will be denoted simply by $A_{i}$. 
Note that because of the ordering of the $A_{i}$, the face $\ext(A_{i})$ is in 
the complex $\Delta_{i-1}$, so each of these subdivisions is defined.

We claim that $\Delta_{k} = \Delta(L\setminus \{\hat{0}\})$.  The proof
proceeds by induction on $n$.  The case $n=1$ is clear.

When $n>1$, consider the complex $\Delta_{1}= sd_{\ext(A_{1})}(\Delta_{0})$.
Since $\Delta_{0}$ is a simplex, the new vertex $A_{1}$ is a cone point, that is, it is
in every maximal simplex of $\Delta_{1}$.  The base of this cone (the {\em link} of $A_{1}$) consists of
all the facets $F_{1}, \dots , F_{m}$ of $\Delta_{0}$ that are opposite to vertices in
$\ext(A_{1})$.  By relabeling if necessary, we can assume that
$\ext(A_{1}) = \{1,2,\dots,m\}$, and $F_{i}= \{ \closure{j}~ |~ j\ne i \}$.  Since all
further subdivisions are made on faces not containing $A_{1}$, the vertex $A_{1}$
remains a cone point in all $\Delta_{i}$.  So it is enough to consider the effect of further
subdivisions on each of the facets $F_{i}$.

Now, by induction, the face $F_{i}$ is subdivided so that it becomes the order complex
of  $L_{i}\setminus \{\hat{0}\}$, where $L_{i}$ is the lattice of closed subsets of
$[n]\setminus \{i\}$.  Since $A_{1}$ is a cone point
in $\Delta_{k}$, and $A_{1}$ is in every maximal chain in $L$, it follows that $\Delta_{k}$
is the order complex of $L\setminus \{\hat{0}\}$.
\end{proof}

\begin{figure}[htb]
\begin{center}
\setlength{\unitlength}{1600sp}%
\begingroup\makeatletter\ifx\SetFigFont\undefined%
\gdef\SetFigFont#1#2#3#4#5{%
  \reset@font\fontsize{#1}{#2pt}%
  \fontfamily{#3}\fontseries{#4}\fontshape{#5}%
  \selectfont}%
\fi\endgroup%
\begin{picture}(13597,2895)(1951,-2161)
\put(13500,-2500){\makebox(0,0)[lb]{\smash{{\SetFigFont{11}{13.2}{\rmdefault}{\mddefault}{\itdefault}$\Delta_3  \cdots  \Delta_6$}}}}
\thicklines
{\color[rgb]{0,0,0}\put(7150,600){\line(-3,-5){1210}}
\put(5950,-1400){\line( 1, 0){2400}}
\put(8350,-1400){\line(-3, 5){1190}}
}%
{\color[rgb]{0,0,0}\put(7150,614){\circle{80}}
}%
{\color[rgb]{0,0,0}\put(5950,-1400){\circle{80}}
}%
{\color[rgb]{0,0,0}\put(8350,-1400){\circle{80}}
}%
{\color[rgb]{0,0,0}\put(11440,-460){\circle{80}}
}%
{\color[rgb]{0,0,0}\put(13750,-460){\circle{80}}
}%
{\color[rgb]{0,0,0}\put(7800,-460){\circle{80}}
}%
{\color[rgb]{0,0,0}\put(15020,-460){\circle{80}}
}%
{\color[rgb]{0,0,0}\put(10150,-470){\circle{80}}
}%
{\color[rgb]{0,0,0}\put(3610,614){\circle{80}}
}%
{\color[rgb]{0,0,0}\put(10800,614){\circle{80}}
}%
{\color[rgb]{0,0,0}\put(14400,600){\circle{80}}
}%
{\color[rgb]{0,0,0}\put(4840,-1400){\circle{80}}
}%
{\color[rgb]{0,0,0}\put(2400,-1400){\circle{80}}
}%
{\color[rgb]{0,0,0}\put(15600,-1400){\circle{80}}
}%
{\color[rgb]{0,0,0}\put(13200,-1400){\circle{80}}
}%
{\color[rgb]{0,0,0}\put(9600,-1400){\circle{80}}
}%
{\color[rgb]{0,0,0}\put(14450,-1400){\circle{80}}
}%
{\color[rgb]{0,0,0}\put(12000,-1400){\circle{80}}
}%
{\color[rgb]{0,0,0}\put(10800,614){\line(-3,-5){1210}}
\put(9600,-1400){\line( 1, 0){2400}}
\put(12000,-1400){\line(-3, 5){1190}}
}%
\put(11400,-490){\line(-2,-1){1820}}
\put(10150,-470){\line( 1, 0){1250}}
{\color[rgb]{0,0,0}\put(14400,600){\line(-3,-5){1210}}
\put(13200,-1400){\line( 1, 0){2380}}
\put(15600,-1400){\line(-3, 5){1190}}
}%
\put(15000,-490){\line(-2,-1){1830}}
\put(13750,-480){\line( 1, 0){1250}}
\put(7800,-470){\line(-2,-1){1830}}
{\color[rgb]{0,0,0}\put(14450,-1400){\line( 3, 5){570}}
}%
\put(3451,-2500){\makebox(0,0)[lb]{\smash{{\SetFigFont{11}{13.2}{\rmdefault}{\mddefault}{\itdefault}$\Delta_0$}}}}
\put(3800,539){\makebox(0,0)[lb]{\smash{{\SetFigFont{11}{13.2}{\rmdefault}{\mddefault}{\itdefault}$\closure{c}$}}}}
\put(7400,539){\makebox(0,0)[lb]{\smash{{\SetFigFont{11}{13.2}{\rmdefault}{\mddefault}{\itdefault}$\closure{c}$}}}}
\put(11000,539){\makebox(0,0)[lb]{\smash{{\SetFigFont{11}{13.2}{\rmdefault}{\mddefault}{\itdefault}$\closure{c}$}}}}
\put(14550,539){\makebox(0,0)[lb]{\smash{{\SetFigFont{11}{13.2}{\rmdefault}{\mddefault}{\itdefault}$\closure{c}$}}}}
\put(1951,-1800){\makebox(0,0)[lb]{\smash{{\SetFigFont{11}{13.2}{\rmdefault}{\mddefault}{\itdefault}$\closure{b}$}}}}
\put(4726,-1800){\makebox(0,0)[lb]{\smash{{\SetFigFont{11}{13.2}{\rmdefault}{\mddefault}{\itdefault}$\closure{a}$}}}}
\put(5776,-1800){\makebox(0,0)[lb]{\smash{{\SetFigFont{11}{13.2}{\rmdefault}{\mddefault}{\itdefault}$\closure{b}$}}}}
\put(8176,-1800){\makebox(0,0)[lb]{\smash{{\SetFigFont{11}{13.2}{\rmdefault}{\mddefault}{\itdefault}$\closure{a}$}}}}
\put(9451,-1800){\makebox(0,0)[lb]{\smash{{\SetFigFont{11}{13.2}{\rmdefault}{\mddefault}{\itdefault}$\closure{b}$}}}}
\put(11776,-1800){\makebox(0,0)[lb]{\smash{{\SetFigFont{11}{13.2}{\rmdefault}{\mddefault}{\itdefault}$\closure{a}$}}}}
\put(13126,-1800){\makebox(0,0)[lb]{\smash{{\SetFigFont{11}{13.2}{\rmdefault}{\mddefault}{\itdefault}$\closure{b}$}}}}
\put(15376,-1800){\makebox(0,0)[lb]{\smash{{\SetFigFont{11}{13.2}{\rmdefault}{\mddefault}{\itdefault}$\closure{a}$}}}}
\put(8026,-436){\makebox(0,0)[lb]{\smash{{\SetFigFont{11}{13.2}{\rmdefault}{\mddefault}{\itdefault}$\closure{a,c}$}}}}
\put(11551,-436){\makebox(0,0)[lb]{\smash{{\SetFigFont{11}{13.2}{\rmdefault}{\mddefault}{\itdefault}$\closure{a,c}$}}}}
\put(15151,-436){\makebox(0,0)[lb]{\smash{{\SetFigFont{11}{13.2}{\rmdefault}{\mddefault}{\itdefault}$\closure{a,c}$}}}}
\put(9200,-436){\makebox(0,0)[lb]{\smash{{\SetFigFont{11}{13.2}{\rmdefault}{\mddefault}{\itdefault}$\closure{b,c}$}}}}
\put(12800,-436){\makebox(0,0)[lb]{\smash{{\SetFigFont{11}{13.2}{\rmdefault}{\mddefault}{\itdefault}$\closure{b,c}$}}}}
\put(14026,-1800){\makebox(0,0)[lb]{\smash{{\SetFigFont{11}{13.2}{\rmdefault}{\mddefault}{\itdefault}$\closure{a,b}$}}}}
\put(7051,-2500){\makebox(0,0)[lb]{\smash{{\SetFigFont{11}{13.2}{\rmdefault}{\mddefault}{\itdefault}$\Delta_1$}}}}
\put(10651,-2500){\makebox(0,0)[lb]{\smash{{\SetFigFont{11}{13.2}{\rmdefault}{\mddefault}{\itdefault}$\Delta_2$}}}}
{\color[rgb]{0,0,0}\put(3600,600){\line(-3,-5){1200}}
\put(2400,-1400){\line( 1, 0){2400}}
\put(4840,-1400){\line(-3, 5){1200}}
}%
\end{picture}%
\caption{The sequence of $\Delta_{i}$ for the example in Figure \ref{threepoints}.\label{deltas}}
\end{center}
\end{figure}
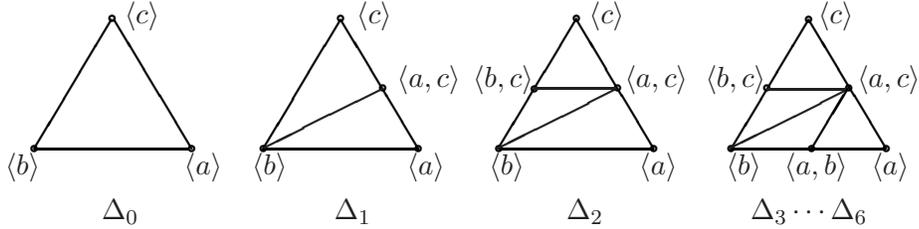

Notice that the stellar subdivisions over the principal closed sets (join-irreducibles) are
redundant and can be omitted without loss. 
Figure \ref{deltas} gives the
sequence of subdivisions leading to the order complex of the meet-distributive lattice
generated by the example of three collinear points.
In this example, once $\Delta_3$ is constructed, every subsequent subdivision is
over a principal closed set and therefore has no effect on the complex.

By a {\em polyhedral ball} we will mean a simplicial complex that
is topologically a $d$-dimensional ball and can be embedded to
give a regular triangulation, that is, one that admits a
strictly convex piecewise-linear function (see, for example,
\cite{BFS} for the definitions).  Polyhedral balls are known to satisfy
strong enumerative conditions \cite{BL}.

\begin{corollary}
For any meet-distributive lattice $L\ne B_n$, the order complex 
$\Delta(L\setminus \{\hat{0},\hat{1}\})$ is a polyhedral
ball.
\end{corollary}
\begin{proof}
Let $L$ be the lattice of closed sets of the convex closure $\closure{\cdot}$
on $[n]$.  Since $L \ne B_n$, we have that $\ext([n])\ne [n]$, and so
every stellar subdivision that is involved in producing 
$\Delta(L\setminus \{\hat{0}\})$ takes place on the boundary of
the simplex $\Delta_0$.

Since being the boundary complex of a simplicial convex polytope is
preserved under taking stellar subdivisions \cite{PrB}, we conclude
that the boundary of $\Delta(L\setminus \{\hat{0}\})$ is the boundary
of a simplicial convex polytope $Q$.  By means of a projective transformation
that sends the vertex $A_1=[n]=\hat{1}$ to the point at infinity, we see
that the image of $Q$ under such a map is the graph of a strictly
convex function over $\Delta(L\setminus \{\hat{0},\hat{1}\})$.
\end{proof}

It was shown in \cite{PrB} that stellar subdivision preserves the property
of being {\em vertex decomposable}, which in turn implies shellability.
As a consequence we get that both $\Delta(L\setminus \{\hat{0}\})$
and $\Delta(L\setminus \{\hat{0},\hat{1}\})$ are vertex
decomposable and hence shellable, as stated in \cite[Theorem 8.1]{BKL}
(and its proof) in the language of greedoids. 
Theorem \ref{ordercomplex}
was first proved for distributive lattices in \cite{Pr} precisely
to show that order complexes of distributive lattices were shellable.
The result in \cite{Pr} was stated for $\Delta(L)$, which is a cone
over the complex we consider.


\section{The associated CW spheres}
We define now a triangulated sphere derived from the order complex of
a meet-distributive lattice $L$.  It will turn out that this triangulated
sphere is the barycentric subdivision of a regular CW sphere that has
the same enumerative relationship to $L^{*}$ (the dual to $L$) as an arrangement of hyperplanes
(oriented matroid) has to the underlying geometric lattice.

\subsection{The complex $\pm\Delta$}

For a meet-distributive lattice $L$, let
$\Delta = \Delta(L\setminus \{\hat{0}\})$, a triangulation of the
$(n-1)$-simplex $\Delta_0$.  We will define a
triangulation $\pm\Delta$ of the $n$-dimensional crosspolytope $O_n$ as
follows.  If the vertex set of the simplex $\Delta_0$ is $[n]=\{1,2,\dots,n\}$,
then that of the crosspolytope is $\pm[n]=\{\pm 1,\pm 2,\dots,\pm n\}$.
Faces of the crosspolytope are all $\sigma \subseteq \pm[n]$ such that
not both $i$ and $-i$ are in $\sigma$.

	We reflect the triangulation $\Delta$ to obtain a triangulation
$\pm\Delta$ of the crosspolytope, much as the crosspolytope can be built
by reflecting the simplex generated by the unit vectors.
We consider the simplex $\Delta_0$ to be embedded as the convex hull of
the unit vectors and define the triangulation $\pm\Delta$ by reflecting
the triangulation $\Delta$.

Formally, $\pm\Delta$ is the simplicial
complex whose vertices are all equivalence classes of
pairs $(A,\varepsilon)$, where
$A\in L\setminus \{\hat{0}\}$, $\varepsilon$ is a map from
$[n]$ to $\{\pm 1 \}$, and we identify $(A,\varepsilon)$ and
$(A,\varepsilon')$ when $\varepsilon|_{\ext(A)}=\varepsilon'|_{\ext(A)}$.
For arbitrary $\varepsilon : [n]\rightarrow \{\pm 1\}$
and $\sigma= \{A_1,A_2,\dots,A_k\}\in \Delta$, let
\begin{equation*}
\sigma^{\varepsilon}:= \big\{(A_1,\varepsilon),
(A_2,\varepsilon), \dots , (A_k,\varepsilon)\big\}
\end{equation*}
and
$$\Delta^{\varepsilon}= \big\{ \sigma^{\varepsilon} | \sigma \in \Delta\big\}.$$
$\Delta^{\varepsilon}$ is essentially the triangulation $\Delta$ transferred
to the face of the crosspolytope given by the sign pattern $\varepsilon$.
Finally, define
$$\pm\Delta = \bigcup_{\varepsilon} \Delta^{\varepsilon},$$
the union being taken over all $\varepsilon : [n]\rightarrow \{\pm 1\}$.

\begin{remark}
Note that boundary faces of $\Delta$ can result in faces of $\pm\Delta$
having more
than one name; in fact, $\sigma^\varepsilon = \sigma^\rho$ if and only if
$\varepsilon$ and $\rho$ agree on the set
$$\ext(\sigma):= \bigcup_{i=1}^k \ext(A_i).$$
\label{facename}
\end{remark}

\begin{figure}[htb]
\begin{center}
\setlength{\unitlength}{2000sp}%
\begingroup\makeatletter\ifx\SetFigFont\undefined%
\gdef\SetFigFont#1#2#3#4#5{%
  \reset@font\fontsize{#1}{#2pt}%
  \fontfamily{#3}\fontseries{#4}\fontshape{#5}%
  \selectfont}%
\fi\endgroup%
\begin{picture}(3624,5124)(1789,-4873)
\thicklines
{\color[rgb]{0,0,0}\put(3601,-1261){\line(-2,-1){1800}}
\put(1801,-2161){\line( 6,-5){1800}}
\put(3601,-3661){\line( 6, 5){1800}}
\put(5401,-2161){\line(-2, 1){1800}}
}%
{\color[rgb]{0,0,0}\put(1801,-2161){\line( 3,-1){1800}}
\put(3601,-2761){\line( 3, 1){1800}}
}%
{\color[rgb]{0,0,0}\put(3601,239){\line( 0,-1){3000}}
\put(3601,-2761){\line( 0,-1){2100}}
}%
{\color[rgb]{0,0,0}\put(3601,239){\line(-3,-4){1800}}
\put(1801,-2161){\line( 2,-3){1800}}
\put(3601,-4861){\line( 2, 3){1800}}
\put(5401,-2161){\line(-3, 4){1800}}
}%
{\color[rgb]{0,0,0}\put(3601,-1261){\line(-3,-4){900}}
\put(2701,-2461){\line( 3,-4){900}}
\put(3601,-3661){\line( 3, 4){900}}
\put(4501,-2461){\line(-3, 4){900}}
}%
{\color[rgb]{0,0,0}\put(2701,-961){\line( 3,-1){900}}
\put(3601,-1261){\line( 3, 1){900}}
}%
{\color[rgb]{0,0,0}\put(3001,-3961){\line( 2, 1){600}}
\put(3601,-3661){\line( 2,-1){600}}
}%
\end{picture}%
\caption{Triangulation of the boundary of the octahedron induced by
reflecting $\Delta_{6}$.}\label{reflect}
\end{center}
\end{figure}
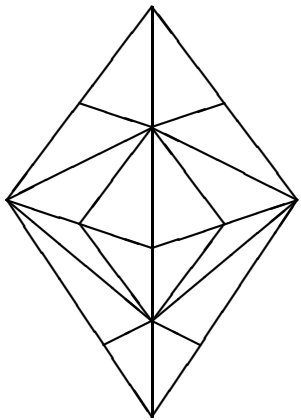

Figure \ref{reflect} shows the complex $\pm\Delta$ for the example of three collinear points.

\begin{theorem}
For any meet-distributive lattice $L$,  $\pm\Delta(L\setminus \{\hat{0}\})$
can be obtained from the $n$-dimensional crosspolytope by a sequence
of stellar subdivisions, and so it is combinatorially the boundary
complex of an $n$-dimensional simplicial polytope, where $n$ is the number of
join-irreducibles of $L$.
\label{spherethm}
\end{theorem}
\begin{proof}
As before, suppose $L$ is the lattice of closed sets of a convex
closure $\closure{\cdot}$ on $[n]$, and let $A_{1}, A_{2}, \dots , A_{k}$
be a reverse linear extension of $L\setminus \{\hat{0}\}$.  
We can extend this order to an order on all pairs $(A_i,\varepsilon)$,
$\varepsilon : [n] \rightarrow \{\pm 1\}$,
by ordering lexicographically,
with the given order on the first coordinate and any order on the second.

It is now relatively straightforward to adapt the proof of Theorem \ref{ordercomplex}
to show that the complex $\pm\Delta$
is obtained by carrying out stellar subdivisions over faces of $O_n$
in the order given by the order of the $(A_i,\varepsilon)$.  The
subdivision corresponding to $(A_i,\varepsilon)$ is done over the
face $\{ \varepsilon(j) \cdot {j}~ |~ j\in \ext(A_i) \}$ of $O_n$; in $\pm\Delta$,
every face containing $\{ \varepsilon(j) \cdot {j}~ |~ j\in \ext(A_i) \}$ is
subdivided as it would be by doing the subdivision in the boundary
of $O_n$.  Again, since stellar subdivision preserves the property of being
the boundary complex of a polytope, the result follows.
\end{proof}

\subsection{The poset $Q_L$} \label{S:CW}

	We construct a regular CW complex $\Sigma$ having $\pm\Delta$
as its barycentric subdivision.  Equivalently, if $\F(\Sigma)$ is the
face poset of $\Sigma$, then $\Delta(\F(\Sigma)) = \pm\Delta$.

We begin by defining a poset $Q_L$ associated to any meet-distributive
lattice $L$.  The elements of $Q_L$ are all equivalence classes of
pairs $(A,\varepsilon)$, where
$A\in L$ and $\varepsilon$ is a map from $[n]$ to $\{\pm 1 \}$ as before.
We define the order relation on $Q_L$ by $(A,\varepsilon) \le (B,\delta)$
if and only if $A \subseteq B$ and the maps $\varepsilon,\delta$ agree
on the set $\ext(A) \cap \ext(B)$.  We include an element $\hat{1}\in Q_L$
for convenience; the element $(\hat{0},\emptyset)$ corresponding to
$\hat{0}\in L$ serves as $\hat{0}$ in $Q_L$.
Note that when $L$ is distributive, $Q_L$ is the signed Birkhoff poset of \cite{H}.

\begin{proposition}
$\Delta(Q_L\setminus \{\hat{0},\hat{1}\}) = \pm\Delta$
\label{pmDelta}
\end{proposition}
\begin{proof}
The maximal simplices in $\pm\Delta$ are the simplices
\begin{equation*}
\sigma^{\varepsilon}:= \{(A_1,\varepsilon),
(A_2,\varepsilon), \dots , (A_n,\varepsilon)\},
\end{equation*}
where
$A_1 \subseteq A_2 \subseteq\cdots\subseteq A_n$
is a maximal chain in $L \setminus \{\hat{0}\}$.  Then clearly,
$$(A_1,\varepsilon)<
(A_2,\varepsilon) < \dots < (A_n,\varepsilon)$$
is a maximal chain in $Q_L \setminus \{\hat{0}, \hat{1}\}$.

Conversely, if
$$(A_1,\varepsilon_1)<
(A_2,\varepsilon_2) < \dots < (A_n,\varepsilon_n)$$
is a maximal chain in $Q_L \setminus \{\hat{0},\hat{1}\}$, then, if we let
$\sigma = \{ A_1,A_2,\dots,A_n \} \in \Delta$, we have $\ext(\sigma)=[n]$
and so there is an $\varepsilon : [n] \rightarrow \{\pm 1\}$ such that
$\varepsilon_i = \varepsilon|_{\ext(A_i)}$ for each $i$.  Thus
$$\{(A_1,\varepsilon_1),
(A_2,\varepsilon_2), \dots , (A_n,\varepsilon_n)\} = \sigma^\varepsilon$$
is a maximal simplex in $\pm\Delta$.
\end{proof}

Next, we define a cell complex $\Sigma_L$ from 
the lattice $L$ (the underlying convex closure $\closure{\cdot}$
on $[n]$) and the simplicial complex $\pm\Delta$
as follows.  For each $A \in L\setminus \{\hat{0}\}$ and 
$\varepsilon :[n] \rightarrow \{\pm 1\}$, we define a
cell $C_{(A,\varepsilon)}$ that is a union of simplices in $\pm\Delta$.  
For $A = [n]$, we take $C_{(A,\varepsilon)}$ to be the star of
$(A,\varepsilon)$ in the complex $\pm\Delta$, that is, the union of
all maximal simplices containing the vertex $(A,\varepsilon)$.

For proper closed sets $A\in L$, we consider the subgeometry
$\closure{\cdot}$ restricted to subsets of
$A$, with lattice $L_A = [\hat{0},A]$ and order complex
$\Delta_A = \Delta(L_A \setminus  \{\hat{0}\} )$.
The complex $\Delta_A$ is the subcomplex of $\Delta$ subdividing the
the face of $\Delta_0$ spanned by the vertices $i\in A$,
and the corresponding complex
$\pm\Delta_A$ is the subcomplex of $\pm\Delta$
subdividing the faces of the crosspolytope spanned by all vertices
$\pm i, i\in A$.  For any $A\in L$ and
any $\varepsilon : [n] \rightarrow \{\pm 1\}$, we define the
cell $C_{(A,\varepsilon)}$ to be the star of $(A,\varepsilon)$
in the complex $\pm\Delta_A$.

Since $\pm\Delta_{A}$ is the boundary of a simplicial polytope by
Theorem \ref{spherethm}, each cell
$C_{(A,\varepsilon)}$ is topologically a disk of dimension $|A|-1$,
and its boundary
is the link of the vertex $(A,\varepsilon)$ in the complex $\pm\Delta_A$
and so is a sphere. 
We define $\Sigma_L$ to be the collection of all the cells
$C_{(A,\varepsilon)}$, $A\in L\setminus \{\hat{0}\}$.

\begin{lemma}
The boundary of
$C_{(A,\varepsilon)}$ is the union of all cells
$C_{(B,\delta)}$, where  $B \subseteq A, B \ne A$ and the maps
$\varepsilon,\delta$ agree on $\ext(A) \cap \ext(B)$.
\label{bdylemma}
\end{lemma}
\begin{proof}
By definition, we have
\begin{equation}
C_{(A,\varepsilon)} = \bigcup_{\stackrel{A\in \sigma \in \Delta}
{\gamma |_{\ext(A)}=\varepsilon  |_{\ext(A)} }}
\sigma^{\gamma}.
\label{Ccell}
\end{equation}
Since $\partial C_{(A,\varepsilon)}$ is the link of $(A,\varepsilon)$ in $\pm\Delta_A$, that is,
\[ 
\partial C_{(A,\varepsilon)} =
 \bigcup_{\stackrel{\tau \in \lk_{\Delta_{A}}(A)}{\gamma |_{\ext(A)}=\varepsilon  |_{\ext(A)}}}
\tau^{\gamma},
\]
the statement of the lemma is equivalent to
\begin{equation}
\bigcup_{\stackrel{\tau \in \lk_{\Delta_{A}}(A)}{\gamma |_{\ext(A)}=\varepsilon  |_{\ext(A)} } }
\tau^{\gamma}= 
\bigcup_{\stackrel{B \subsetneq A}
{\delta |_{\ext(A)\cap \ext(B)}=\varepsilon  |_{\ext(A)\cap \ext(B)}} } 
C_{(B,\delta)}.
\label{Cbdy}
\end{equation}
Here the unions are over $\gamma : [n] \rightarrow \{\pm 1\}$ and
$\delta : [n] \rightarrow  \{\pm 1\}$, respectively,
and $ \lk_{\Delta_{A}}(A) = \{ \tau \in \Delta_{A} ~|~ A\notin \tau, \tau \cup \{A\}
 \in \Delta_{A}\}$ is the link of $A$ in $\Delta_{A}$.

To see the equality in (\ref{Cbdy}), note that if $\tau^\gamma$,
$\tau \in \lk_{\Delta_{A}}(A)$, $\gamma |_{\ext(A)}=\varepsilon  |_{\ext(A)}$, appears
on the left side, then $\tau^\gamma \subseteq C_{(B,\gamma)}$,
where $B$ is a maximal element of $\tau$.  Since
$\gamma|_{\ext(A)\cap \ext(B)}=\varepsilon  |_{\ext(A)\cap \ext(B)}$,
the cell $C_{(B,\gamma|_{\ext(B)})}$ appears on the right side.

For the opposite inclusion, suppose $\tau^\gamma$ is a maximal
simplex of $\pm\Delta_{A}$ in $C_{(B,\delta)}$,
where $B \subsetneq A$ and
$\delta |_{\ext(A)\cap \ext(B)}=\varepsilon  |_{\ext(A)\cap \ext(B)}$.
Then $\gamma|_{\ext(B)}= \delta|_{\ext(B)}$ by (\ref{Ccell}), and so
$$\gamma|_{\ext(A)\cap \ext(B)}=\delta |_{\ext(A)\cap \ext(B)}=
\varepsilon  |_{\ext(A)\cap \ext(B)}.$$
Since $i\in \ext(A) \cap B$ implies $i \in \ext(B)$ (otherwise
$i \in \closure{B\setminus \{i\}} \subseteq \closure{A\setminus  \{i\}}$),
we have that the only places where $\gamma|\ext(A)$ and $\varepsilon$
might not agree are outside of $B$.  Since $\ext(\tau) \subseteq B$,
we may, by Remark \ref{facename},
adjust $\gamma$ to $\gamma'$ outside of $B$ so that
$\tau^{\gamma'} = \tau^\gamma$ and $\gamma'|_{\ext(A)} = \varepsilon|_{\ext(A)}$.
Thus $\tau^\gamma$ appears on the left of (\ref{Cbdy}), establishing the
equality.
\end{proof}

We can now prove the main result of this section.
\begin{theorem}
The cells in $\Sigma_L$ form a regular CW sphere, with face poset
$Q_L\setminus \{\hat{0},\hat{1}\}$  and barycentric subdivision $\pm\Delta$.
\end{theorem}
\begin{proof}
Since each $\Delta^{\varepsilon}$ is a cone on
$([n],\varepsilon)$, 
$$|\Sigma_L| = \bigcup_{(A,\varepsilon)\in Q_L\setminus \{\hat{0},\hat{1}\}} 
C_{(A,\varepsilon)} = 
\bigcup_{\varepsilon:[n] \rightarrow \{\pm 1\}}  C_{([n],\varepsilon)}= |\pm \Delta|,$$
so $|\Sigma_L|$ is a sphere by Theorem \ref{spherethm}.

By construction, the only inclusions 
 $C_{(B,\delta)} \subseteq C_{(A,\varepsilon)}$
possible among cells is when
 $C_{(B,\delta)} \subseteq \partial C_{(A,\varepsilon)}$,
so  $C_{(B,\delta)} \subseteq C_{(A,\varepsilon)}$ if and
only if $(B,\delta) \le (A,\varepsilon)$ in
$Q_L \setminus \{\hat{0},\hat{1}\}$
by Lemma \ref{bdylemma}.

To see that $|\Sigma_L|$ is a regular CW sphere, one can assemble
$|\Sigma_L|$ according to a linear extension of the poset
$Q_L \setminus \{\hat{0},\hat{1}\}$.  By Lemma \ref{bdylemma}, all the
boundary faces
of any cell $C_{(A,\varepsilon)}$  will be present when it comes time
to attach it.

Since the poset of inclusions among the faces of $\Sigma_L$ is
$Q_L \setminus \{\hat{0},\hat{1}\}$, it will have $\pm\Delta$ as
barycentric subdivision by Proposition \ref{pmDelta}.
\end{proof}

\subsection{Join-distributive lattices} \label{S:join-dist}
We note briefly that everything in this section works for
{\em join-distributive lattices}, that is lattices $L$ whose
dual $L^*$ (reverse all order relations) is meet-distributive.
Here we have to reverse the roles of $\hat{0}$ and $\hat{1}$.
In particular, both $L$ and $L^*$ have the same order complex,
so we have $\Delta(L\setminus \{\hat{1}\}) = \Delta(L^* \setminus\{\hat{0}\})
=\Delta$, which gives rise to the same simplicial polytope $\pm\Delta$.

For join-distributive $L$, the poset $Q_L = \left(Q_{L^*}\right)^*$,
and so the corresponding spherical complex $\Sigma_L$ is defined by
defining the maximal cells to correspond to the maximal elements
of $Q_L\setminus\{\hat{1}\}$ (the minimal elements of $Q_{L^*}\setminus\{\hat{0}\}$). Here,
the CW sphere $\Sigma_L = \left(\Sigma_{L^*}\right)^*$ is the dual
to $\Sigma_{L^*}$.

\begin{figure}[htb]
\begin{center}
\resizebox{!}{1.5in}{\includegraphics{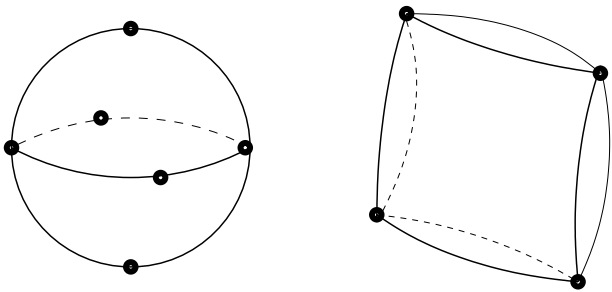}}
\caption{The spheres $\Sigma_{L}$ and $\Sigma_{L^{*}}$ for $L$
from three collinear points.\label{unsubdivide}}
\end{center}
\end{figure}

Figure \ref{unsubdivide} shows the CW sphere for both the meet-distributive $L$
from three collinear points and the corresponding join-distributive $L^{*}$.
Note that both $\pm\Delta$ and $\Sigma_{L}$ retain the full $\left(\Z/2\Z\right)^{n}$
symmetry of the crosspolyope.

We remark that if $L$ is the lattice of the convex geometry on $n$ collinear points,
then one can verify that $Q_{L^*}$ is isomorphic to the Tchebyshev poset 
$T_n$ of \cite{He}.  
Hetyei showed that $T_n$ is the face poset of a
regular CW sphere and its order complex subdivides a crosspolytope. In fact his
proof of the latter assertion uses essentially the reflection construction
discussed at the beginning of \S 3.1.


\section{Enumerative properties of $Q_L$}

For a graded poset $P$ (with $\hat{0}$ and $\hat{1}$) with rank function $\rho$, 
define
\begin{equation} \label{E:nu}
\nu(P) = \sum_{t\in P} (-1)^{\rk(t)}\mu(\hat{0}, t),
\end{equation}
where $\mu$ denotes the M\"obius function as defined in \cite[Chapter 3]{ECI}.
If $L$ is the intersection lattice of a real hyperplane arrangement,
a well known result due to Zaslavsky \cite{Z} gives $\nu(L)$ as the number of 
connected components in the  complement of the arrangement. 
He extended this to show how all the face numbers of an arrangement
depend solely on the lattice of intersections.
As a generalization, \cite[Proposition~4.6.2]{matroid} expresses the flag numbers
of an arrangement, that is, the enumerators of chains of faces having prescribed rank sets, in terms of
the functional $\nu$ applied to intervals in the intersection lattice.

We now show that for  join-distributive $L$,
the flag numbers of $Q_L$ may be computed similarly from
intervals in $L$. (For meet-distributive $L$, the flag numbers of $Q_L$ can be obtained from this
by duality.)  Suppose that
$L$ consists of the closed sets of a convex geometry
ordered by reverse inclusion. In analogy with the zero map on oriented matroids \cite{matroid}, 
we define the map  $z:Q_L \setminus \{\hat{0}\} \to L$ by $z((A,\varepsilon)) = A$. 
 
\begin{proposition}\label{P:birkhoff-zeromap}
Let $c = \{ A_1 <  A_2 < \cdots < A_k = \hat{1}\}$ be a chain in the join distributive lattice
$L$ and $z^{-1}(c)$ denote the set of chains  in $Q_L$ 
that are mapped by $z$ to $c$. Then
\[ |z^{-1}(c)| =  \prod_{i=1}^{k-1} \nu([A_{i}, A_{i+1}]).
\]
\end{proposition}
\begin{proof}
Given a sign function $\varepsilon_i : [n] \to\{\pm 1\}$ for some $2 \le i \le k$, 
there are $2^{|\ext(A_{i-1})\setminus\ext(A_{i})|}$ essentially different sign
functions $\varepsilon_{i-1}:[n] \to \{\pm 1\}$ such that
$(A_{i-1}, \varepsilon_{i-1}) < (A_{i}, \varepsilon_{i})$ in $Q_L$, since the only restriction on
 $\varepsilon_{i-1}$ is that it  agree with $\varepsilon_{i}$ on 
$\ext(A_{i}) \cap \ext(A_{i-1})$. Thus, starting with $\varepsilon_k = \emptyset$,
there are precisely $\prod_{i=1}^{k-1} 2^{|\ext(A_i)\setminus \ext(A_{i+1})|}$ ways to build a sequence of
sign functions $\varepsilon_k, \ldots, \varepsilon_{2},\varepsilon_1$ resulting in a chain
$(A_1, \varepsilon_1) < \cdots < (A_{k},\varepsilon_{k})$ in $Q_L$.

To complete the proof it suffices to show that for $1\le i \le k$,
\begin{equation}\label{E:mobius-sum}
\sum_{A_{i} \le B \le A_{i+1} } (-1)^{\rk(A_i,B)} \mu(A_i,B) = 
2^{|\ext(A_{i}) \setminus \ext(A_{i+1})|}.
\end{equation}
The M\"obius function of  a join-distributive lattice satisfies
\[
(-1)^{\rk(A_i,B)} \mu(A_i,B)=  
\begin{cases}
1 & \text{ if $[A_i,B]$ is a Boolean lattice,} \\
0 & \text{ otherwise.}
\end{cases}
\]
(This follows, for example, from \cite[Theorems~4.2, 4.3]{EdelJam}.)
By definition of join-distributivity, for $B\in [A_{i}, A_{i+1}]$
the interval $[A_i,B]$ is a Boolean lattice precisely when $B$ is less than or equal to 
the join of atoms of $[A_{i}, A_{i+1}]$, which are  those
$A_i \setminus \{a\}$ such that $a \in \ext(A_{i}) \setminus \ext(A_{i+1}).$ Hence
the left side of  \eqref{E:mobius-sum} reduces to a sum of the form $\sum_{B} 1$ with
$B$ ranging over a Boolean lattice of rank  $|\ext(A_{i})\setminus\ext(A_{i+1})|$.
\end{proof}

The complete enumerative information on chains in a graded poset
$P$ is carried by the formal power series
\[
F_P :=  \sum_{\hat{0}=t_0 < t_1 < \cdots  < t_{k}=\hat{1} \atop 0< i_1 < \cdots < i_k}
x_{i_1}^{\rk(t_0,t_1)} x_{i_2} ^{\rk(t_1,t_2)} \cdots x_{i_k}^{\rk(t_{k-1},t_k)},
\]
where  $\rk(s,t) = \rho(t)-\rho(s)$.
As $P$ ranges over the family of graded posets, the
$F_P$ span the (Hopf) algebra of
quasisymmetric functions, denoted $\QSym$. The definition of $F_P$ is due to 
Ehrenborg \cite{Ehr}. 
See \cite[\S 7.19]{ECII} for further background on quasisymmetric functions.

In the context of combinatorial Hopf algebras  \cite{ABS}, the functional
$\nu$ can be seen as the pullback of a certain ``odd character" $\nu_\QSym$ to the 
Hopf algebra of graded posets 
along the map $P \mapsto F_P$; that is, $\nu(P) = \nu_\QSym(F_P)$.
 By  the general theory there is an induced
Hopf algebra map $\vartheta: \QSym \to \QSym$ satisfying
\[
\vartheta(F_P) =   \sum_{\hat{0}=t_0 < t_1 < \cdots  < t_{k}=\hat{1} \atop 0< i_1 < \cdots < i_k}
\nu([t_0, t_1]) \cdots \nu([t_{k-1}, t_k])
x_{i_1}^{\rk(t_0,t_1)} \cdots x_{i_k}^{\rk(t_{k-1},t_k)}.
\]
In fact $\vartheta$ is precisely the map introduced by Stembridge \cite{Ste} to
relate the quasisymmetric weight enumerator for $P$-partitions of a labeled poset to the
enriched quasisymmetric weight enumerator of that poset.
See \cite[Examples~2.2, 4.4, 4.9]{ABS}.

The main result of this section is an extension of \cite[Theorem~5.15]{H}:

\begin{theorem}\label{T:theta-birkhoff}
For a join-distributive lattice $L$, we have
\[ 2 F_{Q_L} =  \vartheta (F_{L \, \cup \, \{\hat{0}\}}),\]
where $\hat{0}$ denotes a new minimum element adjoined to $L$.
\end{theorem}
\begin{proof}
This is essentially the argument used to prove \cite[Theorem~3.1]{BER}.
Extend $z$ to a map $z: Q_L\to
L\cup \{\hat{0}\}$ by requiring  $z(\hat{0}) = \hat{0}$.
By Proposition~\ref{P:birkhoff-zeromap},
\begin{eqnarray*}
F_{Q_L} 
& = & \sum_{c = \{ \hat{0} = A_0  < \cdots < A_k = \hat{1}\} \subseteq L \, \cup \,\{\hat{0}\}
\atop i_1 < \cdots < i_k} |z^{-1}(c)| 
x_{i_1}^{\rk(A_0, A_1)} \cdots x_{i_k}^{\rk(A_{k-1}, A_k)} \\
& = & \sum_{ \hat{0} = A_0 < \cdots < A_k = \hat{1} \atop i_1 < \cdots < i_k} 
\nu([A_1, A_2]) \cdots \nu([A_{k-1},A_k])
x_{i_1}^{\rk(A_0,A_1)} \cdots x_{i_k}^{\rho(A_{k-1},A_k)} \\
& = & \frac{1}{2} \sum_{ \hat{0} = A_0 < \cdots < A_k = \hat{1} \atop i_1 < \cdots < i_k} 
\nu([A_0, A_1]) \cdots \nu([A_{k-1},A_k])
x_{i_1}^{\rk(A_0,A_1)} \cdots x_{i_k}^{\rho(A_{k-1},A_k)} \\
& = & \frac{1}{2} \vartheta(F_{L \, \cup \, \{\hat{0}\}}).
\end{eqnarray*}
The third equality holds because $\hat{0}$ is covered by only one element in $L\cup \{\hat{0}\}$,
implying that $\mu(\hat{0}, A_1)$ vanishes if $\rk(A_1) >1$; hence $\nu([\hat{0}, A_1]) =
 \mu(\hat{0}, \hat{0}) - \mu(\hat{0},[n]) = 2$.
\end{proof}

Analogously, if $Z$ is the face lattice of the zonotope associated with a hyperplane arrangement
and  $L$ is the intersection lattice of the arrangement, then 
\[ 2 F_{Z} =  \vartheta(F_{L \, \cup \, \{\hat{0}\}})\]
\cite{BER} \cite[Proposition~3.5]{BHW}. 
It is easy to see that a join-distributive lattice $L$ must be
semimodular, as are geometric lattices.  
One is led to speculate whether this relationship holds for {\em all} semimodular
lattices, namely, whether 
for any semimodular lattice $L$,  
there exists a regular CW sphere $\Sigma_{L}$
with face poset $Q_L$ (with $\hat{0}$, $\hat{1}$ adjoined) such that
\[ 2 F_{Q_L} = \vartheta( F_{L \, \cup \, \{\hat{0}\}}). \]
The role played by convex closures
in this work might be played instead by {\em interval greedoids}
(see \cite[Theorem 8.8.7]{BZ}; we are grateful to Anders Bj\"orner for
suggesting this connection).  Note that this would imply  the existence of spheres
$\Sigma_{L}$ for geometric lattices that are not necessarily orientable.
In nonorientable case, one might also ask for the relationship of
$\vartheta( F_{L \, \cup \, \{\hat{0}\}})$
to the (dual) face counts of the homotopy-sphere arrangements of Swartz \cite{Sw}.
Simple examples suggest the former might provide lower bounds for the latter.
In the orientable case, these bounds are clearly achieved by the results of \cite{BHW,BER}.
One could speculate further that achieving the bounds implies orientability.

\bigskip

There is a well-known bijection between multichains of a fixed length in a distributive lattice $J(P)$ and $P$-partitions (order-preserving maps) whose parts have a certain fixed upper bound \cite[Proposition~3.5.1]{ECI}. Edelman  and Jamison \cite[Theorem~4.7]{EdelJam} extended this to a bijection between multichains in a meet-distributive lattice and extremal functions, which are generalizations of $P$-partitions. We will conclude with a discussion of an analogous correspondence between multichains in $Q_L$ and a new class of functions called {\it enriched extremal functions,} which are generalizations of enriched $P$-partitions \cite{Ste}. (We are grateful to Paul Edelman for suggesting that we seek such a correspondence.)

Let $L$ be the lattice of closed sets of a convex closure $\closure{\cdot}$ on the set $[n]$. Consider the linear ordering $-1 \prec 1 \prec -2 \prec 2 \prec \cdots$ of the set of non-zero integers $\Z \setminus \{0\}$. For a function $f:[n]\to \Z \setminus \{0\}$ and a closed set $A$, let  $f_A$ denote the minimum element of $\{ f(a)~|~a\in A\}$ with respect to $\prec$. 
\begin{definition}
Given a convex closure on $[n]$, a function $f : [n] \to \Z \setminus \{0\}$ is called {\it enriched extremal function}  provided that
\begin{enumerate}
\item for every closed set $A$ there exists $a \in\ext(A)$ such that $f(a) = f_A$, and
\item for every $a\in [n]$, if $f(a) <0$ then $a\in\ext \{ b \in [n]~|~ f(b) \succeq f(a) \}$.
\end{enumerate}
\end{definition}
We remark that if $f$ satisfies condition (1) then  $\{ a \in [n]~|~f(a) \succeq b\}$ is closed for any $b\in \Z\setminus \{0\}$. This justifies the notation used in (2).

For the convex closure on three collinear points there are four enriched extremal functions $f: \{a,b,c\} \to \{ -1, 1\}$. They are given by $(f(a), f(b), f(c)) = (1,1,1), (-1,1,1), (1,1,-1),$ and $(-1,1,-1)$. 

Notice that if $\closure{\cdot}$ is the upper-order-ideal closure on a poset $P = ([n], \le_P)$, then $f$ is an enriched extremal function if and only if for all $a <_P b$, we have (1) $f(a) \preceq f(b)$ and (2) $f(a) = f(b)$ implies $f(a) > 0$; in other words $f$ is an {\it enriched $P$-partition} with respect to a natural labeling of $P$ \cite{Ste}.

The zeta polynomial of a graded poset $Q$, denoted $Z(Q,t)$, is determined by the property that for a positive integer $m$,  $Z(Q,m)$ is the number of multichains in $Q_L$ of the form $\hat{0} = q_0 \le q_1 \le \cdots \le q_m = \hat{1}$. It will be convenient to introduce another polynomial $\overline{Z}(Q,t)$ given by
\[\overline{Z}(Q,t) = \sum_{q \text{ maximal in } Q \setminus \{\hat{1}\}} Z([\hat{0}, q], t).\]
Thus $\overline{Z}(Q,m)$ is the number of multichains in $Q$ of the form $\hat{0} = q_0 \le \cdots \le q_m$, where $q_m$ is a maximal element in $Q \setminus \{\hat{1}\}$.

For a positive integer $m$, recall $\pm [m] = \{-m, \ldots, -2, -1, 1, 2, \ldots, m\}.$

\begin{proposition} \label{P:enriched} Suppose that  $L$ is the meet-distributive lattice of closed sets of a convex closure $\closure{\cdot}$ on $[n]$. Then for $m \ge 1$,
\[ \overline{Z}(Q_L, m) = \text{\# of enriched extremal functions $f: [n] \to \pm [m]$}.\]
 \end{proposition}
\begin{proof}[Sketch of proof]
Given a multichain $([n], \varepsilon_0) = (A_0, \varepsilon_0) \ge (A_1, \varepsilon_1) \ge \cdots \ge  (A_m , \varepsilon_m) = (\emptyset, \emptyset)$ in $Q_L\setminus \{\hat{1}\}$, we obtain an enriched extremal function $f:[n]\to \{ \pm 1, \ldots, \pm m\}$ by setting, for $a \in A_{i-1} \setminus A_i$,
\[
f(a) = \begin{cases} -i & \text{if $a \in \ext(A_{i-1})$ and $\varepsilon_{i-1}(a) = -1$,} \\ 
i & \text{otherwise.} \end{cases}
\]
Conversely, for an enriched extremal function $f: [n] \to \pm [m]$, if we write $\{ |f(1)|, |f(2)|, \ldots, |f(n)|\} = \{s_1 < s_2 < \cdots < s_k\}$ then we can recover the corresponding multichain by setting $A_i = \{ a \in [n] ~|~ f(a) \succeq  - s_j\}$ for $1\le j \le k$ and $s_{j-1} \le i < s_j$, where $s_0 = 0$, and setting $\varepsilon_{i}(a)$ equal to the sign of $f(a)$.
\end{proof}

To obtain a similar formula for the zeta polynomial $Z(Q_L,m)$, we extend $\closure{\cdot}$ to a convex closure on $[n+1]$ by declaring that $\closure{n+1} = [n+1]$. If $L'$ denotes the lattice of closed sets for the new closure then it is easy to see that $\overline{Z}(Q_{L'},m) = 2 \; Z(Q_L, m)$. It follows that
\[ 2 \;Z(Q_L, m) = \text{\it \# of enriched extremal functions $f: [n+1] \to \pm [m]$}.\]

Proposition~\ref{P:enriched} may be viewed as the enriched analogue of \cite[Theorem~4.7]{EdelJam}, which asserts that $Z(L,m)$ enumerates certain extremal functions and that by reciprocity $(-1)^n Z(L,-m)$ enumerates {\it strictly} extremal functions. It turns out that for $Q_L$ we have self-reciprocity, that is, 
\[Z(Q_L, - t) = (-1)^{n+1} Z(Q_L, t)\]
and
\[\overline{Z}(Q_L, - t) = (-1)^{n} \overline{Z}(Q_L, t)\]
for any meet-distributive lattice $L$ (cf.\ \cite[Proposition~4.2]{Ste}). This follows, for example, from \cite[Proposition~3.14.1]{ECI} and the fact that $Q_L$ is an Eulerian poset of rank $n+1$.

A further corollary of Proposition~\ref{P:enriched} is that if $L$ is the lattice of upper order ideals of a poset $P$ and $P_0$ denotes the poset obtained from $P$ by adjoining a new minimum element, then 
\begin{equation}\label{E:enriched} \nonumber
 2 \; Z(Q_L, m) = \text{\it \# of enriched $P_0$-partition $f:P_0 \to \pm [m]$.}
 \end{equation}
(See also \cite[Corollary~5.3]{H}.) As an application we describe a way to translate the recent counterexamples of Stembridge's enriched poset conjecture \cite{Ste04} into new counterexamples of Gal's real root conjecture for flag triangulations of spheres \cite{G}.

It can be shown (e.g., \cite[Chapter~3, Exercise~67b]{ECI}) that the generating function for $Z(Q_L, m)$ satisfies
\[ \sum_{m \ge 0} Z(Q_L, m)\;  t^m = 
 \frac{t \cdot h^{\Delta(Q_L\setminus\{\hat{0},\hat{1}\})}(t)}{(1-t)^{n+1}}, \]
where $h^{\Delta(Q_L\setminus\{\hat{0},\hat{1}\})}(t)$ denotes the $h$-{\it polynomial} of the order complex $\Delta(Q_L\setminus\{\hat{0},\hat{1}\})$. Stembridge found examples of a poset $P$ such that the numerator of the rational generating function enumerating enriched $P_0$-partitions has non-real roots, thereby disproving an earlier conjecture of his \cite{Ste}. From such a poset one can construct via our results a flag simplicial
complex (meaning every minimal non-face has size two) -- namely
$\Delta(Q_L \setminus\{\hat{0}, \hat{1}\})$, where $L$ is the lattice of upper order ideals of $P$ -- that barycentrically subdivides
a regular CW sphere and whose $h$-polynomial has non-real roots.
In fact, this simplicial sphere will be a simplicial polytope. We have thus provided additional counterexamples of Gal's conjecture that the $h$-polynomial of a flag simplicial triangulation of a sphere should have only real roots  \cite{G}.

\providecommand{\bysame}{\leavevmode\hbox
to3em{\hrulefill}\thinspace}

\end{document}